\documentclass[12pt]{elsarticle}
\usepackage[T2A]{fontenc}
\usepackage[cp1251]{inputenc}
\usepackage[english]{babel}
\usepackage{amsmath,amssymb,amsthm}

\newtheorem{theorem}{Theorem}

\newdefinition{remark}{Remark}
\newdefinition{example}{Example}

\def\FF{{\mathcal F}}
\def\RR{{\mathbb R}}
\def\uu{{\bf u}}
\def\vv{{\bf v}}
\def\ww{{\bf w}}
\def\0{{\bf 0}}
\def\xx{{\bf x}}
\def\ee{{\bf e}}
\let\emp\varnothing

\begin{document}

\title{Finding a subset of nonnegative vectors with a coordinatewise large sum\tnoteref{note}}
\tnotetext[note]{Supported by the Russian government project 11.G34.31.0053.}

\author[mipt,yar]{Ilya I.~Bogdanov\corref{cor}}
\ead{ilya.i.bogdanov@gmail.com}
\cortext[cor]{Corresponding author}

\author[yar]{Grigory R.~Chelnokov}
\ead{grishabenruven@yandex.ru}

\address[mipt]{Moscow Institute of Physics and Technology (State University), Institutsky per., 9,
Dolgoprudny, Moscow reg., Russia 141700}

\address[yar]{Laboratory of Discrete and Computational Geometry, Yaroslavl'
State University, Sovetskaya st. 14, Yaroslavl', Russia 150000}

\begin{abstract}
  Given a rational $a=p/q$ and $N$ nonnegative $d$-dimensional real vectors ${\bf u}_1$, \dots, ${\bf u}_N$, we show that it is always possible to choose $(d-1)+\left\lceil   (pN-d+1)/q\right\rceil$ of them such that their sum is (componentwise) at least  $(p/q)({\bf u}_1+\dots+{\bf u}_N)$. For fixed $d$ and $a$, this bound is sharp if $N$ is large enough. The method of the proof uses Carath\'eodory's theorem from linear programming.
\end{abstract}


\maketitle

\section{Introduction}

We deal with the $d$-dimensional real vector space~$\RR^d$; the vectors of the standard basis
are denoted by $\ee_1,\dots,\ee_d$. Introduce a coordinatewise
partial order~$\succeq$ on $\RR^d$; that is, for the vectors
$\uu=[u^1,\dots,u^d]$ and $\vv=[v^1,\dots,v^d]$ we write $\uu\succeq
\vv$ if $u^j\geq v^j$ for $1\leq j\leq d$.

Let $\uu_1,\dots, \uu_N\in \RR^d$ be $N$ nonnegative vectors
(that is, $\uu_i\succeq \0$ for $1\leq i\leq N$), and let $a\in[0,1]$ be some
real number. We say that a set of indices $I\subseteq \{1,\dots,N\}$ is
{\em $a$-rich} if
$$
  \sum_{i\in I}\uu_i\succeq a\sum_{i=1}^N \uu_i.
$$

Let $f_{N,d}(a)$ be the minimal number $f$ such that for every $N$ nonnegative vectors
$\uu_1,\dots,\uu_N\in \RR^d$ there exists an $a$-rich set $I$ with $|I|\leq f$.
Further we consider only rational $a$ and write $a=p/q$ with $q> 0$ and $\gcd(p,q)=1$.

To find an upper bound for $f_{N,d}(a)$, one may use a theorem of Stromquist and Woodall~\cite{str-wood} claiming that, given $n$ non-atomic probability measures on~$S^1$, there exists a union of $n-1$ arcs that has measure $a$ in each measure. It can be performed as follows. Let $w^j=\sum_{i=1}^N u_i^j$; we may assume that $w_j>0$ for $1\leq j\leq d$. Consider a segment $T=[0,N]$, identify its endpoints to obtain a circle of length~$N$, and split it into unit segments. For $1\leq i\leq N$ and $1\leq j\leq d$, define a measure $\mu^j$ on segment $[i-1,i]$ as $\mu^j=u_i^j\mu/w^j$, where $\mu$ is the usual Lebesgue measure; set also $\mu^{d+1}=\mu/N$. By the theorem mentioned above, there exists a union of $d$ arcs $J\subseteq T$ such that $\mu^j(\FF)=a$ for $1\leq j\leq d+1$. Now, one may define
$$
  I=\left\{i:J\cap [i-1,i]\neq \emp\right\}.
$$
This set is $a$-rich since
$$
  \sum_{i\in I}u_i^j=w^j\mu^j\left(\bigcup\nolimits_{i\in I}[i-1,i]\right)\geq
  w^j\mu^j(J)=a\sum_{i=1}^n u_i^j.
$$
Moreover,
$$
  |I|\leq \mu(J)+2d=N\mu^{d+1}(J)+2d= aN+2d.
$$
Thus, $f_{N,d}(a)\leq aN+2d$.

In an analogous way, one may apply a well-known Alon's
theorem on splitting of necklaces~\cite{alon} obtaining a bound
$$
  f_{N,d}(p/q)\leq \frac pq\cdot N+\frac{p(q-p)}q\cdot d.
$$

The bounds shown above are asymptotically tight. Nevertheless, they provide
exact values of~$f_{N,d}(p/q)$ only for some border cases. The aim of this paper
is to find an exact value of~$f_{N,d}(a)$ for every rational $a=p/q$, positive
integer $d$ and sufficiently large $N$. We use only the methods of linear
programming.

The main result is the following theorem.

\begin{theorem}
  For any positive integer numbers $N$, $d$ and rational number $a=p/q\in[0,1]$,
  we have
  $$
    f_{N,d}(p/q) \le (d-1)+\left\lceil\frac {pN-d+1}q\right\rceil.
  $$

  Moreover, if $q>p\geq 1$ and $N\geq (q-1)(d-1)$, then we have
  $$
    f_{N,d}(p/q) = (d-1)+\left\lceil\frac {pN-d+1}q\right\rceil.
  $$
%
  \label{main}
\end{theorem}

Throughout the rest of the paper, we use the notation $s=d-1$.

\smallskip
The next section contains the proof of the upper bound. Here we present an
example showing that this bound is sharp if $N$ is large enough.

\begin{example} Choose an integer $r\in[1,q-1]$ such
that $pr\equiv 1\pmod q$. Let $m=\lceil pr/q\rceil$ (hence $qm-pr=q-1$).

Let us set $\uu_{ir-k}=\ee_i$ for $1\leq i\leq s$ and
$0\leq k\leq r-1$, and set $\uu_i=\ee_d$ for $i>rs$ (notice that $N\geq (q-1)s\geq rs$). Denote
$\ww=\sum_{i=1}^N\uu_i=[r,r,\dots,r,N-rs]$. Now, if $\sum_{i\in I}
\uu_i\succeq \frac pq\ww$, then
$$
  \sum_{i\in I}\uu_i\succeq \left[m,m,\dots,m,\left\lceil \frac
  pq(N-rs)\right\rceil\right],
$$
because all the coordinates of $\uu_i$ are integer. Thus, since the sum of
coordinates of each vector is~1, we should have
$$
  |I|\geq ms+\left\lceil \frac pq(N-rs)\right\rceil
  =s+\left\lceil \frac {pN+s(qm-pr-q)}q\right\rceil
  =s+\left\lceil \frac {pN-s}q\right\rceil,
$$
as desired.
\qed
\end{example}

\section{Proof of the upper bound}

Consider $N$ vectors $\uu_1,\dots,\uu_N\in\RR^d$ with nonnegative coordinates. Denote
$$
  \ww=\sum_{i=1}^N\uu_i,
  \qquad
  f=s+\left\lceil\frac {pN-s}q\right\rceil.
$$
We need to prove that there exists a set of indices $I$ such that
$$
  |I|\leq f
  \qquad \text{and} \qquad
  \sum_{i\in I}\uu_i\succeq \frac pq\ww.
$$
We use induction on $p+q$. In the base cases $p+q\leq 2$ we have $a=0$ or $a=1$, and the statement is trivial.

Now, assume that $p+q\geq 3$, and assume that the statement of the theorem holds for all
pairs $(p',q')$ with $p'+q'<p+q$. Introduce the following set:
$$
  X_{p/q}=\left\{\xx\in \RR^N\colon\quad 0\leq x_i\leq 1,\quad \sum_{i=1}^N
  x_i\uu_i=\frac pq \ww\right\}.
$$
This set is closed and bounded; next, it is nonempty since $[p/q,\dots,p/q]\in X_{p/q}$.
Moreover, this set is defined by $s+1$ linear equalities and some linear
inequalities. By Carath\'eodory's theorem, there exists a
vector $\xx=[x_1,\dots,x_N]\in X_{p/q}$ such that $N-s-1$ of these inequalities
come to equalities; that is, $N-s-1$ coordinates of~$\xx$ are integer. Hence,
either $(i)$ at least $N-f$ coordinates are zeros, or $(ii)$ at least
$f-s$ coordinates are ones.

\smallskip
In case $(i)$, denote $I=\{i: x_i>0\}$. We have $|I|\leq N-(N-f)=f$. On the
other hand, we obtain
\begin{equation}
  \sum_{i\in I}\uu_i\succeq \sum_{i\in I} x_i\uu_i=\sum_{i=1}^N x_i\uu_i=\frac
  pq\ww,
  \label{sum-i}
\end{equation}
as desired.
\smallskip

In case $(ii)$, define $J=\{i: x_i <1\}$, and let $N'=|J|$. We have
\begin{equation}
  \sum_{i\in J}\uu_i=\ww-\sum_{i\notin J}\uu_i
  \succeq \ww-\sum_{i=1}^N x_i\uu_i=\frac{q-p}q\ww.
  \label{sum-j}
\end{equation}
Notice that in~\eqref{sum-i} and~\eqref{sum-j} we have used the condition that all the vectors~$\uu_i$ are nonnegative.

Next, note that
\begin{equation}
  N'=|J|\leq N-f+s=N-\left\lceil\frac{pN-s}q\right\rceil\leq
  N-\frac{pN-s}q=\frac{(q-p)N+s}q.
  \label{N'}
\end{equation}
Renumbering the vectors we may assume that $J=\{1,2,\dots,N'\}$. Again, we
distinguish two subcases: $(ii')$ $q\geq 2p$ and $(ii'')$ $q<2p$.

\smallskip
In case $(ii')$, we apply the induction hypothesis to the vectors
$\uu_1,\dots,\uu_{N'}$ and the number $a'=p/(q-p)\in[0,1]$. We obtain the
subset $I\subseteq \{1,\dots,N'\}$ such that
$$
  |I|\leq s+\left\lceil\frac {pN'-s}{q-p}\right\rceil
$$
and
$$
  \sum_{i \in I} \uu_i
  \succeq \frac{p}{q-p}\sum_{i \in J} \uu_i.
$$
By~\eqref{sum-j}, the last inequality yields
\begin{equation}
  \sum_{i \in I} \uu_i\succeq
  \frac{p}{q-p}\cdot \frac{q-p}{q}\ww=\frac{p}{q}\ww.
  \label{ii'-1}
\end{equation}
Next,  by~\eqref{N'} we have
$$
  pN'-s\leq \frac pq\bigl((q-p)N+s\bigr)-s=\frac{q-p}q(pN-s),
$$
thus obtaining
\begin{equation}
  |I| \leq s+\left\lceil\frac 1{q-p}\cdot\frac{q-p}q(pN-s)\right\rceil=f.
  \label{ii'-2}
\end{equation}
The relations~\eqref{ii'-1} and~\eqref{ii'-2} show that $I$ is a desired set
of indices.

\smallskip
In case $(ii'')$, we apply the induction hypothesis to $N-N'$ vectors
$\uu_{N'+1}$, \dots, $\uu_{N}$ and the number $a'=(2p-q)/p\in(0,1)$. We obtain
the subset $I'\subseteq \{N'+1,\dots,N\}$ such that
$$
  |I'|\leq s+\left\lceil\frac {(2p-q)(N-N')-s}p\right\rceil
$$
and
$$
  \sum_{i \in I'} \uu_i
  \succeq \frac{2p-q}p\sum_{i=N'+1}^{N} \uu_i.
$$

Now we claim that the subset $I=I'\cup J$ satisfies the desired properties.
Recall that by~\eqref{sum-j} we have
$$
  \sum_{i=1}^{N'}\uu_i=\frac{q-p}q\ww+\ww'
$$
for some vector $\ww'\succeq \0$. Hence
$$
  \sum_{i=N'+1}^{N}\uu_i=\frac pq\ww-\ww',
$$
and we obtain
\begin{multline*}
  \sum_{i\in I}\uu_i=\sum_{i\in I'}\uu_i+\sum_{i\in J}\uu_i
  \succeq
  \frac{2p-q}p\left(\frac pq\ww-\ww'\right)+\left(\frac{q-p}q\ww+\ww'\right)\\
  =\frac pq\ww+\frac{q-p}p\ww'\succeq\frac pq\ww.
\end{multline*}
We are left to show that $|I|\leq f$.

Recall that
\begin{multline*}
  |I|=|J|+|I'|\leq N'+s+\left\lceil\frac {(2p-q)(N-N')-s}p\right\rceil\\
  =s+\left\lceil\frac{(q-p)N'+(2p-q)N-s}p\right\rceil.
\end{multline*}
So, it suffices to prove that
$$
  \frac{(q-p)N'+(2p-q)N-s}p\leq \frac{pN-s}q,
  \qquad\text{or} \qquad
  qN'\leq (q-p)N+s,
$$
which is equivalent to~\eqref{N'}. Thus $I$ is a desired set.
\qed

\section*{Acknowledgements}
The authors are grateful to Prof. V.L.~Dolnikov and Prof. R.N.~Karasev for useful discussions. The authors thank the referees for their valuable remarks.

%

\end{document}